\documentclass{elsart}
\usepackage{amssymb}
\usepackage{amsmath}
\usepackage{amsfonts}
\usepackage{graphicx}%
\theoremstyle{definition}
\newtheorem{definition}{Definition}[section]

\newtheorem{rema}[definition]{Remark}
\newtheorem{exa}[definition]{Example}
\newtheorem{obs}[definition]{Observation}

\theoremstyle{plain}

\newtheorem{theorem}[definition]{Theorem}
\newtheorem{corollary}[definition]{Corollary}

\newenvironment{remark}%
{\begin{rema}\rm}%
{\end{rema} }

\newenvironment{example}%
{\begin{exa}\rm}%
{\end{exa} }

\newenvironment{observation}%
{\begin{obs}\rm}%
{\end{obs} }

\newenvironment{proof}{\noindent{\bf Proof.}}{\hfill$\Box$}

\newenvironment{ac}{\noindent{\bf Acknowledgements:}}{}

\def\sg{\mbox{sg}}
\def\Res{{\rm{Res}}}

\def\sres{\mbox{sres}_k}
\def\d{{\bf d}}
\def\x{{\bf x}}
\def\fh{{f}^h}
\def\tfh{\tilde{f}^h}
\def\Proj{{\mathbb P}}

\def\RR{{\mathcal R}}
\def\SS{{\mathcal S}}
\def\TT{{\mathcal T}}
\def\HH{{\mathcal H}}
\def\EE{{\mathcal E}}
\def\VV{{\mathcal V}}

\def\OO{{\mathcal O}}

\def\C{{\mathbb C}}

\def\N{{\mathbb N}}
\def\Q{{\mathbb Q}}

\def\Z{{\mathbb Z}}

\begin{document}

\begin{frontmatter}

\title{Multivariate Subresultants in Roots}

\author{Carlos D' Andrea\thanksref{00}}
\address{Department of Mathematics, University of California at Berkeley, CA 94720 USA.}
\thanks[00]{Research supported by the Miller Institute
for Basic Research in Science, UC Berkeley}
\ead{cdandrea@math.berkeley.edu}

\author{Teresa Krick\thanksref{11}}
\address{Departamento de Matem\'atica, Facultad de
Ciencias Exactas y Naturales, Universidad de Buenos Aires,  1428
Buenos Aires, Argentina}
\thanks[11]{Research supported by grants CONICET PIP
2461/01 and UBACYT X-112.}
\ead{krick@dm.uba.ar}

\author{Agnes Szanto\corauthref{21}\thanksref{22}}

\address{Department of Mathematics, North Carolina State
University, Raleigh, NC 27695 USA}
\thanks[22]{ Research supported by NSF
grants CCR-0306406 and CCR-0347506.}
\corauth[21]{Corresponding author.}
\ead{aszanto@ncsu.edu}



\begin{abstract}
We give a rational expression for the subresultants of
$n+1$ generic polynomials $f_1, \ldots, f_{n+1}$  in $n$ variables
as a function of  the coordinates of the common roots of $f_1,
\ldots, f_n$ and their evaluation in $f_{n+1}$. We present a
simple technique to prove our results, giving new proofs and
generalizing the classical Poisson product formula for the
projective resultant, as well as the expressions of Hong for
univariate subresultants in roots.
\end{abstract}

\end{frontmatter}


\vspace{-7mm}
\section{Introduction}
\vspace{-7mm}
The classical Poisson product formula for resultants of univariate
polynomials can be stated as follows: if $ f$ and $g$ are  two
univariate polynomials of degrees $d_1$ and $d_2$ respectively, with
$g=b_{d_2}(x-\xi_1)\cdots (x-\xi_{d_2})$, then the resultant of $f$ and
$g$ can be expressed as
\begin{equation}\label{poisson}
\Res(f,g)=(-1)^{d_1 d_2}\, b_{d_2}^{d_1}\prod_{j=1}^{d_2}f(\xi_j).
\end{equation}
The main result of this paper is a generalization of Formula
(\ref{poisson}) for univariate  and multivariate subresultants (see
Theorems \ref{subresultant1} and  \ref{subresultant2}). Although most
of the results in the univariate case already appeared in
\cite{Hong99a,Hong01,Las03,DTGV04}, here we present  simple
techniques that enable us to reobtain them    (see Theorem
\ref{subresultant1} and Corollary \ref{hong}) and allow us to
 generalize them to the multivariate case.

Resultants and subresultants of two univariate polynomials go back
to Leibniz, Euler, B\'ezout and Jacobi.  We refer to \cite{GaGe99}
for  historical references. In their modern form, subresultants
were introduced by Sylvester in \cite{Syl}. They have been used to
give an efficient and parallelizable algorithm for computing the
greatest common divisor of two polynomials
\cite{Coll,BrTr,GLRR,duc,HW,HY,rei}. More recently they were also
applied in symbolic-numeric computation
\cite{EmGaLo97,Zeng2003,KalMay2003,ZengDay2004}.

Multivariate resultants were mainly introduced by Macaulay in
\cite{Mac2}, after earlier work by Euler, Sylvester and Cayley,
while multivariate subresultants were first defined by Gonzalez-Vega
in \cite{Gon2,Gon}, generalizing Habicht's method  \cite{Hab}. The
notion of subresultants that we use in the present paper was
introduced by Chardin in \cite{Cha}. It works as follows: let
$f^h_1,\dots,f^h_s$ be a system of generic homogeneous polynomials
in $K[x_0,x_1,\dots,x_n]$ of degrees $d_i=\deg(f^h_i)$ with
parametric coefficients, where $s\leq n+1$ and $K$ is the
coefficient field of $f^h_1,\dots,f^h_s$. Let
$\HH_{d_1,\dots,d_s}:\N\to\N$ be the Hilbert function of a complete
intersection given by $s$ homogeneous polynomials in $n+1$ variables
of degrees $d_1,\dots,d_s.$ Fix $t\in\N$
 and let $\SS$ be a set of $\HH_{d_1,\dots,d_s}(t)$ monomials of degree
$t.$ The \textit{subresultant} $\Delta_\SS$ is a polynomial in $K$
whose degree in the coefficients of $f^h_i$ is
$\HH_{d_1,\dots,d_{i-1},d_{i+1},\dots,d_s}(t-d_i)$ for $i=1,\ldots, s,$
 having the following universal property: $\Delta_\SS$ vanishes at a
particular coefficient specialization $\tfh_1,...,\tfh_{s}\in \C[x_0,
\ldots, x_{n}]$  if and only if $I_t\cup  \SS$ does not generate the
space of all forms of degree $t$. Here, $I_t$ is the degree $t$ part of
the ideal generated by the $\tfh_i$'s (see \cite{Cha}).

 The constructions in \cite{Gon2,Cha}
generalize the classical univariate subresultants in the sense that
they provide the coefficients of certain polynomials in $I_t$, which
in the univariate case include  the greatest common divisor of two
given polynomials.

Theoretical properties and applications of multivariate
subresultants are active areas of research. A series of recent
publications explored: their application to solve zero dimensional
\cite{Gon} and over-constrained polynomial systems \cite{Sza2},  in
the inverse parametrization problem of rational surfaces
\cite{BuDa20042}; their irreducibility and connection with residual
resultants \cite{BuDa2004}; the generalization of their universal
properties to the affine well-constrained case \cite{DaJe2004}; as
well as generalizations of matrix constructions for subresultants
\cite{Sza4}.

Multivariate subresultants also encapsulate as a particular case
the classical projective resultant $\Res(\fh_1, \ldots,
\fh_{n+1})$, which is defined to be an irreducible polynomial in
the coefficients of the $\fh_i$'s which vanishes at a particular
coefficient specialization $\tfh_1,...,\tfh_{n+1}\in \C[x_0,
\ldots, x_{n}]$ if and only if $\tfh_1,...,\tfh_{n+1}$ have a
common root in the complex projective space $\Proj^{ n}_{\C}$.

There is an affine interpretation of the resultant  that can be
stated as follows: Set $$\begin{array}{ccc}
f_i:=\fh_i(1,x_1,\ldots,x_n), & \overline
f_i:=\fh_i(0,x_1,\ldots,x_n), &i=1,\ldots n+1.
\end{array}
$$
Due to B\'ezout's Theorem, the cardinality of the set
$$V(f_1,\ldots,f_n):=\{\xi\in\overline K^n:
f_1(\xi)=f_2(\xi)=\ldots=f_n(\xi)=0\}$$  equals
 $d_1\ldots d_n$ (here, overline denotes algebraic closure),
and the classical Poisson product formula \cite{Weber1912,CoLiSh2,jou},
which generalizes (\ref{poisson}), states that the following identity holds in $\overline{K}$
\begin{equation}\label{poismul}
\Res(\fh_1,...,\fh_{n+1})={\Res(\overline f_1,\ldots,\overline f_n)}^{d_{n+1}}
\prod_{\xi\in V(f_1,\ldots,f_n)}f_{n+1}(\xi).
\end{equation}
In order to make this formula a generalization of (\ref{poisson}),
we have to define resultants of non-homogeneous polynomials. The
obvious generalization is
$\Res(g_1,\ldots,g_{n+1}):=\Res(g^h_1,...,g^h_{n+1})$, where
$g^h_j$ is the homogenization of $g_j$. The same extension to
affine polynomials holds for subresultants.  It should  also be
mentioned that the Poisson formula (\ref{poismul}) is  a
particular case of the determinant of a multiplication map in a
quotient ring (see \cite[Prop. $2.7$]{jou}).
\par
 In Theorem \ref{subresultant2} we generalize (\ref{poismul}) and give
an   expression for \textit{any} multivariate subresultant as a
ratio of two determinants times a function of the coefficients of
$\overline f_1, \ldots, \overline f_n$. The determinant in the
denominator is a Vandermonde type determinant depending on the
common roots of $f_1, \ldots, f_{n}$, while the determinant in the
numerator depends on evaluations of the common roots of $f_1,
\ldots, f_{n}$ in the last polynomial $f_{n+1}$.

\par The paper is structured as follows:   in Section \ref{univariate},
we present in detail the univariate case, showing how to derive with
our techniques Hong's expressions for subresultants of two
univariate polynomials in the roots of one of them and the
coefficients of the other. The details in the univariate case are
essential for the generalization to the multivariate case: they
allow to identify the extraneous factor which is non-trivial in the
multivariate case and they also allow to handle the generality of
the monomial sets appearing in the definition of multivariate
subresultants. In Section \ref{multivariate}, we deal with the
general case.

In order  to keep coherence with the classical literature and previous
works, the presentation in  the univariate case is done in the
traditional way, i.e. for non-homogeneous polynomials, while in the
multivariate case the reader should be aware that the notions involve
homogeneous polynomials.

\begin{ac}
This work started during a visit of the first two authors to North
Carolina State University in the Fall Semester $2004$. We would
like to thank the Department of Mathematics for their hospitality
and for the stimulating working atmosphere it provided. A previous
version of this paper was presented in Mega'05. We thank  the
anonymous referees of that previous version for introducing us to
discrete Wr\'onskians   as presented in \cite{Las03} and for
pointing us other  related references.
\end{ac}

\section{The univariate case}\label{univariate}

{\bf Classical scalar and polynomial subresultants}

 We review here the definition and some well-known properties of the
classical univariate resultant and scalar subresultants and polynomial subresultants.

 Let $f=a_{d_1}x^{d_1}+ \ldots+a_{0} $ and
$g:=b_{d_2}x^{d_2}+\ldots+b_{0}=b_{d_2}(x-\xi_1)\cdots
(x-\xi_{d_2})$, $a_{d_1}\ne 0$, $b_{d_2}\ne 0$, be two polynomials
of degrees $d_1$ and $d_2$ respectively with coefficients in a
field $K$ and roots in the algebraic closure $\overline K$.

The \emph{scalar subresultant} $S_{k}^{(j)}$ of $f$ and $g$ is
defined for $0\leq j\leq j\leq\min\{d_1,d_2\}$ as the following
determinant:

\begin{equation}\label{scalar}
S_{k}^{\left(  j\right)  }:=\det%
\begin{tabular}
[c]{cc}%
$d_1+d_2-2k$ & \\\cline{1-1}%
\multicolumn{1}{|c}{$%
\begin{array}
[c]{cccccc}%
a_{d_1} & \cdots &  & \cdots & a_{k+1-\left(  d_2-k-1\right)  } & a_{j-(d_2-k-1)}\\
& \ddots &  &  & \vdots & \vdots\\
&  & a_{d_1} & \cdots & a_{k+1} & a_{j}%
\end{array}
$} & \multicolumn{1}{|c}{$d_2-k$}\\\cline{1-1}%
\multicolumn{1}{|c}{$%
\begin{array}
[c]{cccccc}%
b_{d_2} & \cdots &  & \cdots & b_{k+1-(d_1-k-1)} & b_{j-(d_1-k-1)}\\
& \ddots &  &  & \vdots & \vdots\\
&  & b_{d_2} & \cdots & b_{k+1} & b_{j}%
\end{array}
$} & \multicolumn{1}{|c}{$d_1-k$}\\\cline{1-1}%
\end{tabular}
\end{equation}
where $a_{\ell}=b_{\ell}=0$ for $\ell<0$.

 The \emph{subresultant
polynomial} $\operatorname*{Sres}\nolimits_{k}(f,g)$ is defined
for $0\le k\le \min\{d_1,d_2\}$ as
$$
\operatorname*{Sres}\nolimits_{k}(f,g):=\sum_{j=0}^{k}S_{k}^{(j)}x^{j}.
$$

 When $k=0$, $\operatorname*{Sres}\nolimits_{0}(f,g)=S_0^{(0)}$ coincides with the classical
 {\em resultant} $\Res(f,g)$
 which arose historically  when checking if $f$ and $g$
have a common factor:
$$\gcd(f,g)= 1 \iff \Res(f,g)\ne 0. $$

In an analogous way, the scalar subresultants satisfy the
following property:
$$\deg\,  \gcd(f,g)=k \ \iff \ S_\ell^{(\ell)}=0 \mbox{ for } 0\le \ell< k
\mbox{  and } S_k^{(k)}\ne 0,$$ and the  polynomial subresultants
$\sres(f,g)$ are  determinant expressions for modified remainders
in the Euclidean algorithm. In particular, for the first $k$ such
that $S_k^{(k)}\ne 0$,  the monic $\gcd$ of $f$ and $g$ satisfies:
$$\gcd(f,g)=(S_k^{(k)})^{-1}\sres(f,g).$$

There is a generalization of the univariate Poisson product formula
(\ref{poisson}) for the polynomial subresultant $\sres (f,g)$, as shown
by Hong in \cite[Th.~3.1]{Hong99a}, see also
\cite[Formula~9.3.2]{Las03} and \cite[Sec.~5]{DTGV04}:
\begin{equation}\label{hongf}\sres(f,g)=(-1)^{(d_1-k)(d_2-k)}b_{d_2}^{d_1-k}
\frac{\left|\begin{array}{ccc}(x-\xi_1)\xi_1^0&\cdots&
(x-\xi_{d_2})\xi_{d_2}^0\\ \vdots & & \vdots \\[-2mm]
(x-\xi_1)\xi_1^{k-1}&\cdots&
(x-\xi_{d_2})\xi_{d_2}^{k-1}\\[-2mm]
 \xi_1^0f(\xi_1) & \cdots & \xi_{d_2}^0f(\xi_{d_2})\\[-2mm]
\vdots & & \vdots \\[-2mm]
 \xi_1^{d_2-k-1}f(\xi_1) & \cdots & \xi_{d_2}^{d_2-k-1}f(\xi_{d_2})\end{array}
 \right|}{ \left|\begin{array}{ccc} \xi_1^0 & \dots & \xi_{d_2}^0\\[-2mm]
 \vdots & & \vdots\\[-2mm]
 \xi_1^{d_2-1}& \dots& \xi_{d_2}^{d_2-1} \end{array}\right| }.
\end{equation}
(Here the sign is due to the fact that we consider $f$ on the roots of
$g$ instead of $g$ on the roots of $f$ as done in \cite{Hong99a}.)

\newpage
 {\bf Notations:}

As we mentioned earlier, most of the results we obtain in this
section are not new. However, we consider important to illustrate
our technique by   applying it  to the  univariate
 case,
 since it  helps to understand its generalization  to the
multivariate setting.  The choices of notations we made here are
accordingly motivated by their coherence with the multivariate
case. They correspond to Chardin's notion of subresultants
\cite{Cha} applied to the univariate case, a slight generalization
of the usual notion of scalar subresultants.

\begin{itemize}
 \item
$f:=a_0+a_1x+\ldots+a_{d_1}x^{d_1} $ and
$g:=b_0+b_1x+\ldots+b_{d_2}x^{d_2}$ in $ K[x]$, where
$K:=\Q(a_0,\dots,a_{d_1},b_0,\dots,b_{d_2})$, with
$a_0,\dots,a_{d_1}, b_0,\dots,b_{d_2}$ algebraically independent
variables over $\Q$ (representing the indeterminate coefficients
of two generic polynomials $f$ and $g$ of degrees $d_1$ and $d_2$
respectively).

 \item
$\{\xi_1,\dots, \xi_{d_2}\}$ denotes the set of roots of $g$ in
$\overline K$ (recall that overline denotes algebraic closure), and
$\VV_{d_2}:=\det(\xi_j^{i-1})_{1\le i,j\le d_2}$ the Vandermonde
determinant associated to this set.

 \item For any $j\in\Z,\
K[x]_j:=\{0\} \cup \{f\in K[x]: \deg f\le j\}.$ Note that if $j<
0,$ then $K[x]_j=\{0\}.$

 \item We set  $t\in\Z$ such that $\,0\leq t\leq d_1+d_2-1$, and let
$t^*:=\max\{d_2-1, t\}$.

 \item $M_f\in K^{(t-d_1+1)\times (t^*+1)}$ and $M_g\in K^{(t-d_2+1)\times (t^*+1)}$
denote  the transposes of the matrices in the monomial bases of
the composition of the Sylvester multiplication maps and the
inclusion $K[x]_t\to K[x]_{t^*}$:
$$\begin{array}{cccc}
\mu_f:& K[x]_{t-d_1}& \to &K[x]_{t^*} \\
& x^\alpha&\mapsto&x^\alpha f(x)\end{array}\quad \mbox{and} \quad
\begin{array}{cccc}
\mu_g:& K[x]_{t-d_2}& \to &K[x]_{t^*} \\
& x^\beta&\mapsto&x^\beta g(x)\end{array},$$ where the monomials
indexing the rows and columns of these matrices are ordered
``increasingly'' $1,\,x,\,x^2,\,\ldots.$ Namely

{\small $$\begin{array}{cc}
 M_f  =    \left[ \! \begin{array}{ccccc|c} a_0\! &
\!\dots \!&\! a_{d_1}\! &\!  \! &\!  \! &\!  \!\\[-2mm]
 \! \! &\! \ddots \! & \!\!&\!  \ddots \!&\!\! &\;  \mathbf{0} \\[-2mm]
 \! \! & \! \! &\!  a_0\!  &\!  \dots \!  & \! a_{d_1} &\!\!
 \end{array}\!  \right]  , &
M_g  =  \left[ \! \begin{array}{ccccc|c}  b_0\! &
\!\dots \!&\! b_{d_2}\! &\!  \! &\!  \! &\!  \!\\[-2mm]
 \! \! &\! \ddots \! & \!\!&\!  \ddots \!&\!\! &\; \mathbf{0} \\[-2mm]
 \! \! & \! \! &\!  b_0\!  &\!  \dots \!  & \! b_{d_2} &\!\!
 \end{array}\!  \right].
 \end{array}$$}
Note that if $t< d_1$ then $M_f=\emptyset$ (the empty matrix), and if
$t<d_2$ then $M_g=\emptyset$.

 \item  We set
 \begin{equation}\label{k}
\begin{array}{ccl}
k&:=&t+1-\dim\big(K[x]_{t-d_1}\big)-\dim\big(K[x]_{t-d_2}\big)\\
&=& t+1 - \max\{0,t-d_1+1\} - \max\{0,t-d_2+1\}\\
&=&t+1 - \max\{0,t-d_1+1\} - (t^*-d_2+1).
\end{array}
\end{equation}
Note that $k\ge 0$ since $t\le d_1+d_2-1$.

 \item
$\SS:=\{x^{\gamma_1},\ldots,x^{\gamma_k};\, 0\le \gamma_1<\dots
<\gamma_k\le t\}\subset K[x]_{t},$ a fixed set of $k$ monomials of degree
bounded by $t$.

 \item $\sg (\SS):=(-1)^\sigma$ where $\sigma$ is a
number of transpositions needed to bring $(1,x,x^2,\ldots,x^{t^*})$ to
$$(x^{\gamma_1},\ldots,x^{\gamma_k},x^{t+1},\ldots,x^{t^*},1,x,
\ldots,x^{\gamma_1-1},x^{\gamma_1+1},\ldots,x^{\gamma_2-1},x^{\gamma_2+1},\ldots,x^t).$$

\item $\Delta_{\SS}:=\Delta^{(t)}_{\SS}(f,g)$ denotes   the {\em
order $t$ subresultant of $f,g$ with respect to $\SS$}, i.e. the
determinant of the matrix whose $\max\{0,t-d_1+1\}$ first rows are
$M_f$, whose $\max\{0,t-d_2+1\}$ following rows are $M_g$ and from
which one deletes the $k+t^*-t$ columns indexed by
$\SS\cup\{x^{t+1},\ldots,x^{t^*}\}.$
\end{itemize}

\begin{remark}\label{relation}

The order $t$ subresultant of $f,g$ with respect to $\SS$
coincides (up to a sign) with  the scalar subresultant when making
special choices of $t$ and  $\SS$:

\begin{enumerate}
\item When $t=d_1+d_2-1$, then $k=t+1-d_2-d_1=0$ and
$\SS=\emptyset$. In that case, from the definitions of $\Res(f,g)$
and   $\Delta_\emptyset$ one gets  that $ \Delta_\emptyset=
(-1)^{d_1 d_2}\,\Res(f,g).$

\item  For $0\le k\le \min\{d_1,d_2\}$ and $t:=d_1+d_2-k-1$, we
can take  $\SS_j:=\{x^i, 0\le i\le k, i\ne j\}$. In that case,
from the definition of $\Delta_{\SS_j}$ and (\ref{scalar}) one
gets that
$\Delta_{\SS_j}=
(-1)^{(d_1-k)(d_2-k)}S_k^{(j)}.$
\end{enumerate}
\end{remark}

The main statement  of this section corresponds to (a slight
generalization of) Hong's theorem \cite[Th.~3.1]{Hong01}. It
expresses $\Delta_\SS$ as the ratio of discrete Wr\'onskians:  we
refer to \cite[Sec.~9.3]{Las03} for an introduction to the
subject. Here we present a new  simple proof  of this result, that
we generalize in the next section to the multivariate setting.

\begin{theorem}\label{subresultant1}
Let $f,g\in K[x]$ and $\{\xi_1,\dots, \xi_{d_2}\}$ be the set of
roots of $g$ in $\overline K$. Then, under the previous notations,
for any fixed $t,\, 0\le t\le d_1+d_2-1$, and  for any
 $\SS=\{x^{\gamma_1},\ldots,x^{\gamma_k}\}\subset K[x]_t $ of cardinality $k,$ with $k$
  defined in (\ref{k}), the order $t$ subresultant $\Delta_\SS$ of $f, g$
  with respect to $\SS$ satisfies:
$$\Delta_\SS= \sg(\SS)\, {b_{d_2}^{t^*-d_2+1}}\,
\frac{|\OO_\SS|}{\VV_{d_2}},$$ where
$$\OO_\SS=\left[\begin{array}{ccc}
\xi_1^{\gamma_1} &\cdots&  \xi_{d_2}^{\gamma_1}\\[-2mm]
\vdots & & \vdots\\[-2mm]
\xi_1^{\gamma_k} &\cdots&  \xi_{d_2}^{\gamma_k}\\
& &  \\[-8mm]\hline & & \\[-8mm]
\xi_1^{t+1} &\cdots&  \xi_{d_2}^{t+1}\\[-2mm]
\vdots & & \vdots\\[-2mm]
\xi_1^{t^*} &\cdots&  \xi_{d_2}^{t^*}\\
& &  \\[-8mm]\hline
& & \\[-8mm]
\xi_1^0f(\xi_1)&\cdots& \xi_{d_2}^0f(\xi_{d_2})\\[-2mm]
\vdots& &\vdots\\[-2mm]
 \xi_1^{t-d_1}f(\xi_1)&
\cdots& \xi_{d_2}^{t-d_1}f(\xi_{d_2})
\end{array}\right] \ \in \ \overline{K}^{d_2\times d_2}.$$
\end{theorem}

\begin{proof}
First,  $\OO_\SS$ is a square matrix since by  (\ref{k}) we have
$$d_2=k+(t^*-t) +\max\{0,t-d_1+1\}.$$
Let $I_{\SS}\in K^{(k+t^*-t)\times (t^*+1)} $ be the transpose of the
matrix of the immersion of the $K$-vector space generated by
$\SS\cup\{x^{t+1},\ldots,x^{t^*}\}$ into $K[x]_{t^*}$ ($I_{\SS}$ is an
identity $(k+t^*-t)$-square matrix plugged into $(t^*+1)$ zero
columns), and set
\begin{equation}\label{Emese}
M_\SS:=\left[\begin{array}{ccc}\;&\;\;\; I_\SS \;\;\;&\;\\ \hline &M_f& \\
\hline& M_g& \end{array} \right].\end{equation} Since it is
straightforward to check  by  (\ref{k}) that we have
$$
k+t^*-t+\max\{0,t-d_1+1\} + \max\{0,t-d_2+1\}= t^*+1,
$$
therefore $M_\SS$ is a $(t^*+1)$-square matrix.

Furthermore, it is immediate to verify  that $| M_\SS|=\sg(\SS)
\Delta_{\SS}$, and we are left to prove that
$|M_\SS|=b_{d_2}^{t^*-d_2+1}\,  |\OO_\SS|/\VV_{d_2}$.

We set $$V_{t^*}:=\left[\begin{array}{ccc} \xi_1^0& \cdots & \xi_{d_2}^0\\[-2mm]
\vdots & & \vdots\\[-2mm]
\xi_1^{t^*}& \cdots & \xi_{d_2}^{t^*}\end{array}\right] \ \in \
\overline{K}^{(t^*+1)\times d_2}, \
V_{d_2}:=\left[\begin{array}{c|c}&{\bf 0}
\\[-2mm]\;\;\; V_{t^*}\;\;\;&\\[-2mm]
 &\;\;Id\;\; \end{array}\right]
 \in\overline{K}^{(t^*+1)\times (t^*+1)}$$ and we observe that
$\VV_{d_2}=| V_{d_2}|$. Now, we perform the product $M_\SS
\,V_{d_2}$:
$$ M_\SS\,V_{d_2}\,=\,
\begin{tabular}{|ccc|}
\cline{1-3}
&$I_\SS $&\\
\cline{1-3}
&$M_f$ & \\
\cline{1-3}
&$M_g$&\\
\cline{1-3}
\end{tabular}
\cdot
\begin{tabular}{|c|c|}
\cline{1-2} &${\bf 0}$ \\$\xi_j^{i-1}$
& \\
 &$Id$
\\
\cline{1-2}
\end{tabular}=\begin{tabular}{|c|ccc|}
\cline{1-4}
& & & \\[-8mm]
 $\xi_j^{\gamma_i}$ &  &  $\bf{*}$ & \\
& & & \\[-8mm]
\cline{1-4}
& & & \\[-8mm]
$\xi_j^{t+i}$ &  &  $\bf{*}$ & \\[-8mm]
& & & \\ \cline{1-4}
& && \\[-8mm]
 $\xi_j^{i-1}f(\xi_j)$ & & ${\bf *}$& \\[-8mm]
& && \\ \cline{1-4}
& & & \\[-8mm]
& $b_{d_2} $& & ${\bf 0}$  \\[-3mm]
$ {\bf 0}$& & $\!\!\ddots\!\!$&   \\[-3mm]
 & ${\bf *}$
&&$b_{d_2} $ \\[-8mm]
& & & \\
\cline{1-4}
\end{tabular}\ .
$$

Therefore  $ | M_\SS |\,\VV_{d_2}= b_{d_2}^{t^*-d_2+1}|\OO_\SS|$,
which proves the Theorem.
\end{proof}

\smallskip

The following examples illustrate how the formula works in a
couple of  cases.


\begin{example}
$d_1=5,\,d_2=2, \, t=4.$ Now we have $t=t^*, \, k=2,$ and
$$\begin{array}{ccc}
M_f=\emptyset, && M_g=\left[\begin{array}{ccccc}
b_0&b_1&b_2&0&0 \\[-2mm]
0&b_0&b_1&b_2&0\\[-2mm]
0&0&b_0&b_1&b_2
\end{array}\right].
\end{array}$$
Set $\SS:=\{x,x^4\}$. Here $\Delta_\SS$ does not coincide with any of
the scalar subresultants $S_2^{(j)}$, $0\le j\le 2$. However, it is
straightforward to check that $\Delta_\SS=b_0b_1^2-b_0^2b_2$. On the
other hand, since $\sg(\SS)=1$, by Theorem \ref{subresultant1} we have
that
 \vspace{-7mm}\begin{eqnarray*} \sg(\SS) \, b_2^{4-2+1}\,\frac{
\left|\begin{array}{cc}
\xi_1&\xi_2\\[-2mm]
\xi_1^4&\xi_2^4
\end{array}\right|}{\left|\begin{array}{cc}
1&1\\[-3mm]
\xi_1&\xi_2
\end{array}\right| }&=& b_2^3\,\xi_1\xi_2\frac{\xi_2^3-\xi_1^3}{\xi_2-\xi_1}\\
&=& b_2^3\,\xi_1\xi_2[(\xi_1+\xi_2)^2-\xi_1\xi_2]\\[1mm]
&=& b_2^3\,(b_0/b_2)\, [(b_1/b_2)^2-(b_0/b_2)] .
\end{eqnarray*}
{\hfill\mbox{$\Box$}}
\end{example}

Next example deals with a case when $t<d_2$ in which case we need
to use  $t^*=d_2-1$ instead of $t$.

\begin{example}
$d_1=2,\,d_2=5, \, t=3$. Here  $k=2$. The scalar subresultants
associated to this value of $k$ are $S_2^{(2)}=a_2^3$,
$S_2^{(1)}=a_2^2a_1$ and $S_2^{(0)}=a_2^2a_0$, while for $t=3<d_2$
 we have $t^*=d_2-1=4$. Thus we have
$$
M_f=\left[\begin{array}{ccccc}
a_0&a_1&a_2&0&0 \\[-2mm]
0&a_0&a_1&a_2&0
\end{array}\right], \quad
M_g=\emptyset. $$ For $\SS:=\{1,x\},$    $\Delta_\SS=a_2^2$, and
Theorem \ref{subresultant1} still works in  this case: since
$\sg(\SS)=1$ and $b_5^{4-5+1}=1$, one has
$$
\frac{ \left|\!\begin{array}{ccccc}
\!1\!&1\!&1\!&1\!&1\!\\[-3mm]
\!\xi_1\!&\xi_2\!&\xi_3\!&\xi_4\!&\xi_5\!\\[-2mm]
\!\xi_1^4\!&\xi_2^4\!&\xi_3^4\!&\xi_4^4\!&\xi_5^4\! \\[-2mm]
\!f(\xi_1)\!&f\!(\xi_2)\!&f(\xi_3)\!&f(\xi_4)\!&f(\xi_5)\!\\[-2mm]
\xi_1\!f(\xi_1)&\xi_2\!f(\xi_2)&\xi_3\!f(\xi_3)&\xi_4\!f(\xi_4)&
\xi_5\!f(\xi_5)
\end{array}\!\right|}{\left|\!\begin{array}{ccccc}
1&1&1&1&1\\[-3mm]
\xi_1&\xi_2&\xi_3&\xi_4&\xi_5\\[-2mm]
\xi_1^2&\xi_2^2&\xi_3^2&\xi_4^2&\xi_5^2\\[-2mm]
\xi_1^3&\xi_2^3&\xi_3^3&\xi_4^3&\xi_5^3\\[-2mm]
\xi_1^4&\xi_2^4&\xi_3^4&\xi_4^4&\xi_5^4
\end{array}\!\right|}
= a_2^2\frac{\left|\!\begin{array}{ccccc}
1&1&1&1&1\\[-3mm]
\xi_1&\xi_2&\xi_3&\xi_4&\xi_5\\[-2mm]
\xi_1^4&\xi_2^4&\xi_3^4&\xi_4^4&\xi_5^4 \\[-2mm]
\xi_1^2&\xi_2^2&\xi_3^2&\xi_4^2&\xi_5^2\\[-2mm]
\xi_1^3&\xi_2^3&\xi_3^3&\xi_4^3&\xi_5^3
\end{array}\!\right|}{\left|\!\begin{array}{ccccc}
1&1&1&1&1\\[-3mm]
\xi_1&\xi_2&\xi_3&\xi_4&\xi_5\\[-2mm]
\xi_1^2&\xi_2^2&\xi_3^2&\xi_4^2&\xi_5^2\\[-2mm]
\xi_1^3&\xi_2^3&\xi_3^3&\xi_4^3&\xi_5^3\\[-2mm]
\xi_1^4&\xi_2^4&\xi_3^4&\xi_4^4&\xi_5^4
\end{array}\!\right|}.
$$
{\hfill\mbox{$\Box$} }
\end{example}


We end this section by  showing  how simple it is to derive from
Theorem \ref{subresultant1} both the  Poisson product formula
 (\ref{poisson}) and  Hong's formula (\ref{hongf}) for  subresultant
polynomials in roots,  together with its generalization to a larger
class of
 determinant polynomials that we call here generalized subresultant
 polynomials.

\begin{observation} \ {\em (Poisson product formula)}
Applying the previous theorem to    Remark \ref{relation}(1), one
obtains
 \vspace{-8mm}\begin{eqnarray*}
\Res(f,g)&=&(-1)^{d_1d_2} \Delta_\emptyset \\
&= &(-1)^{d_1d_2} b_{d_2}^{d_1}\,\frac{ \left|\begin{array}{ccc}
\xi_1^0f(\xi_1)&\cdots& \xi_{d_2}^0f(\xi_{d_2})\\[-2mm]
\vdots& &\vdots\\[-2mm]
 \xi_1^{d_2-1}f(\xi_1)&
\cdots& \xi_{d_2}^{d_2-1}f(\xi_{d_2})
\end{array}\right|
}{ \left|\begin{array}{ccc}
\xi_1^{0} &\cdots&  \xi_{d_2}^{0}\\[-2mm]
\vdots & & \vdots\\[-2mm]
\xi_1^{d_2-1} &\cdots&  \xi_{d_2}^{d_2-1}
\end{array}
\right|}\\ & = & (-1)^{d_1d_2}
b_{d_2}^{d_1}\,\prod_{j=1}^{d_2}f(\xi_j).\end{eqnarray*}
{\hfill\mbox{$\Box$} }
\end{observation}

\begin{observation}\label{hong} {\em (\cite[Th.~3.1]{Hong99a},
\cite[Id.~9.3.2]{Las03})} We derive  Hong's Formula (\ref{hongf})
applying Theorem
 \ref{subresultant1} to
Remark \ref{relation}(2):
 \vspace{-8mm}\begin{eqnarray*}\sres(f,g)&=&\sum_{j=0}^k S_k^{(j)}x^j = (-1)^{(d_1-k)(d_2-k)}
\sum_{j=0}^k \Delta_{\SS_j} x^j \\ &=&
 (-1)^{(d_1-k)(d_2-k)}b_{d_2}^{d_1-k}\VV_{d_2}^{-1}\sum_{j=0}^k
 \sg(\SS_j)|\OO_{\SS_j}|\,x^j.\end{eqnarray*}
 We observe that in this case $t^*=t$,
$\sg(\SS_k)=\sg\{1,\dots,x^{k-1}\}=1$ and $\sg(\SS_j)=(-1)^{k-j}$, and
thus, by column expansion of the determinant we get:
 $$
\sum_{j=0}^k \sg(\SS_j)|\OO_{\SS_j}|\,x^j=\left|\begin{array}{cccc}
(-1)^{k}& \xi_1^0& \dots & \xi_{d_2}^{0}\\[-2mm]
(-1)^{k} x&\xi_1^1& \dots & \xi_{d_2}^{1}\\[-2mm]
\vdots& \vdots &  & \vdots\\[-2mm]
(-1)^{k}x^k& \xi_1^k &\dots & \xi_{d_2}^k\\[-2mm]
0&\xi_1^0f(\xi_1)&\cdots& \xi_{d_2}^0f(\xi_{d_2})\\[-2mm]
\vdots&\vdots& &\vdots\\[-2mm]
0& \xi_1^{d_2-k-1}f(\xi_1)& \cdots& \xi_{d_2}^{d_2-k-1}f(\xi_{d_2})
\end{array}\right|$$
$$= \ (-1)^k  \left|\begin{array}{cccc}
1& \xi_1^0& \dots & \xi_{d_2}^{0}\\[-2mm]
0&\xi_1^1-x\xi_1^0& \dots & \xi_{d_2}^{1}-x\xi_{d_2}^0\\[-2mm]
\vdots& \vdots &  & \vdots\\[-2mm]
0& \xi_1^k-x\xi_1^{k-1} &\dots & \xi_{d_2}^k-x\xi_{d_2}^{k-1}\\[-2mm]
0&\xi_1^0f(\xi_1)&\cdots& \xi_{d_2}^0f(\xi_{d_2})\\[-2mm]
\vdots&\vdots& &\vdots\\[-2mm]
0& \xi_1^{d_2-k-1}f(\xi_1)& \cdots& \xi_{d_2}^{d_2-k-1}f(\xi_{d_2})
\end{array}\right|$$
$$= \left|\begin{array}{ccc}
(x-\xi_1)\xi_1^0& \dots & (x-\xi_{d_2})\xi_{d_2}^{0}\\[-2mm]
\vdots &  & \vdots\\[-2mm]
(x-\xi_1) \xi_1^{k-1}&\dots & (x-\xi_{d_2})\xi_{d_2}^{k-1}\\[-2mm]
\xi_1^0f(\xi_1)&\cdots& \xi_{d_2}^0f(\xi_{d_2})\\[-2mm]
\vdots& &\vdots\\ \xi_1^{d_2-k-1}f(\xi_1)& \cdots&
\xi_{d_2}^{d_2-k-1}f(\xi_{d_2})
\end{array}\right|.
$$
{\hfill\mbox{$\Box$}}
\end{observation}

\bigskip
 One can
straightforwardly generalize Hong's result to a larger class of
determinant polynomials \begin{equation}
\label{chardinpol}
s(x):=\sum_{j=0}^k \Delta_{\SS_{j}}x^{\gamma_j},
\end{equation}
 corresponding to an arbitrary set of
monomials $\SS:=\{x^{\gamma_j},0\le j\le k\}\subset K[x]_t$ and
$\SS_{j}:=S \smallsetminus \{ x^{\gamma_j}\}$, where $d_2\le t\le d_1+d_2-1$
and $k:=d_2 - \max\{0,t-d_1+1\}$. We call such a polynomial a {\em
generalized subresultant polynomial}.

The usual proof that shows that $\sres (f,g)$
 belongs to the ideal $(f,g)$ generated by $f$ and
$g$ extends to showing that $s\in (f,g)$ and the following
expression in terms of roots holds (we omit the proof which is
essentially the same than the proof of Observation \ref{hong}).

\begin{corollary}\label{koko} Let $f,g\in K[x]$ and
$s(x)$ be the generalized
subresultant polynomial   defined in  (\ref{chardinpol}). Then,
   we have
{\small $$s(x)\!=\!b_{d_2}^{t-d_2+1}  \VV_{d_2}^{-1}\!
x^{\gamma_0}\!\left|\begin{array}{ccc} \!
(x^{\gamma_1-\gamma_0}-\xi_1^{\gamma_1-\gamma_0})\xi_1^{\gamma_0}\!&\!
\!\!\!\cdots \! \!\!\!&\!
(x^{\gamma_1-\gamma_0}-\xi_{d_2}^{\gamma_1-\gamma_0})\xi_{d_2}^{\gamma_0}\!\!\!
\\[-2mm]
\!\!\!\vdots \!& \! \!&\! \vdots \!\\[-2mm]\!
(x^{\gamma_k-\gamma_{k-1}}-\xi_1^{\gamma_k-\gamma_{k-1}})\xi_1^{\gamma_{k-1}}
\!&\cdots &
(x^{\gamma_k-\gamma_{k-1}}-\xi_{d_2}^{\gamma_k-\gamma_{k-1}})\xi_{d_2}^{\gamma_{k-1}}
\!\\[-2mm]\!
 \xi_1^0f(\xi_1) \!& \!\!\!\!\cdots \! \!\!\!& \!\xi_{d_2}^0f(\xi_{d_2})\!\\[-2mm]\!
\vdots \!&\! \! & \!\vdots\!\!\! \\[-2mm]\! \xi_1^{d_2-k-1}f(\xi_1) \!&
\!\!\!\!\cdots\! \!\!\!&
\!\xi_{d_2}^{d_2-k-1}f(\xi_{d_2})\!\!\!\end{array}
 \right|.
$$}
{\hfill\mbox{$\Box$}}

\end{corollary}

\section{The multivariate case}\label{multivariate}

\
In this section we generalize Theorem \ref{subresultant1} to Chardin's
multivariate subresultants \cite{Cha}, after introducing the notations
we need.

{\bf Notations:}

\begin{itemize}
\item For  $n\in \N$ and  $1\le i\le n+1$,
$$f_i:=\sum_{|\alpha|\le d_i}a_{i\alpha}{\bf x}^\alpha \ \in \ K[\bf
x],$$ where $\alpha=(\alpha_1,\dots,\alpha_n)\in \left(\Z_{\geq0}\right)^n$,
${\bf x}^\alpha:=x_1^{\alpha_1}\cdots x_n^{\alpha_n}$,
$|\alpha|=\alpha_1+\cdots + \alpha_n$, and $K:=\Q(a_{i\alpha},
1\le i\le n+1, |\alpha|\le d_i)$, with $a_{i\alpha}$ algebraically
independent variables over $\Q$ (representing the indeterminate
coefficients of $n+1$ generic polynomials in $n$ variables $f_i$
of degrees $d_i$ respectively).

 \item For any $j\in\Z$, $ K[{\bf
x}]_j:=K[x_1,\dots,x_n]_j=\{0\} \cup \{f\in K[{\bf x}]: \deg f\le j\}.$

\item  We set $t\in\N$,  $\rho:=(d_1-1)+\cdots + (d_n-1)$  and
$t^*:=\max\{\rho,t\}$.
  \item $k:=\HH_{d_1\dots d_{n+1}}(t)$, the
Hilbert function at $t$ of a regular sequence of $n+1$ homogeneous
polynomials in $n+1$ variables of degrees $d_1,\dots,d_{n+1}$, i.e.
$$k:=\#\{\x^\alpha: |\alpha|\le t,
 \alpha_i<d_i, 1\le i\le n,\, \mbox{and}\, t-|\alpha|<d_{n+1}\}.$$
\item $\SS:=\{\x^{\gamma_1},\ldots,\x^{\gamma_k}\}\subset K[\x]_t$
 a set of $k$ monomials of degree bounded by $t$.

\item For $1\le i\le n+1$,
 $$\RR_i:=\{\x^\alpha, |\alpha|\le t-d_i,
\alpha_j<d_j \ \mbox{for } j<i\}.$$ We observe that for $1\le i\le
n$,
$$\#(\RR_i)=\#\{\x^\alpha,|\alpha|\le t, \alpha_j<d_j \ \mbox{for } j<i
\ \mbox{and} \ \alpha_i\ge d_i\},$$
and
$$\#(\RR_{n+1})=\#\{\x^\alpha,|\alpha|\le t, \alpha_j<d_j \ \forall j
\ \mbox{and} \ t-|\alpha|\ge d_{n+1}\}.$$ Therefore
\begin{eqnarray}\label{MSsquare}
N:={t+n\choose n}=\dim_KK[\x]_t=k+\#(\RR_1)+\cdots +
\#(\RR_{n+1}).
\end{eqnarray}

\item In particular, we denote
$\RR_{n+1}=:\{\x^{\beta_1},\dots,\x^{\beta_r}\},$ where
$r:=\#(\RR_{n+1})$ and we observe that
\begin{equation}\label{key}
k+r=\#\{\x^\alpha,|\alpha|\le t, \alpha_j<d_j \ \forall j\}=\dim
K[{\bf x}]_t/(f_1, \ldots, f_n)\cap K[{\bf x}]_t.
\end{equation}

\item For $j\ge 0$, $\tau_j:= \HH_{d_1\ldots d_n}(j)$, the Hilbert
function at $j$ of a regular sequence of $n$ homogeneous polynomials in
$n$ variables of degrees $d_1,\ldots , d_n$, i.e.
$$\tau_j:=\#\{\x^\alpha: \,|\alpha|=j,
\alpha_i<d_i \mbox{ for } 1\le i\le n\}.$$ We note that $\tau_j=0$
if $j>\rho$.

\item For $j\ge 0$,
\begin{equation}\label{prof}
\TT_{j}:=\left\{\small\begin{array}{l} \mbox{\textit{any} set of }
\tau_j \mbox{ monomials of degree } j
 \ \mbox{ for  }  j\ge \max\{0,t-d_{n+1}+1\},
\\[1mm]
\{\x^\alpha: \,|\alpha|=j, \alpha_i<d_i \mbox{ for } 1\le i\le
n\}\  \mbox{ for  } 0\le j<t-d_{n+1}+1.\end{array}\right.
\end{equation}
See Remark \ref{kokk} for a discussion on the definition of $\TT_j.$

\item $ \TT:=\cup_{j\geq0}\TT_j$ and
$\TT^*:=\cup_{j=t+1}^{t^*}\TT_j$.
 We note that $\#\TT={\d}$,
 where ${\d} :=d_1\cdots d_n$ is the {\em B\'ezout number}, the number of common solutions of
$f_1,\ldots,f_n$ in $\overline{K}^n$, and that $\TT^*=\emptyset$
if $t^*=t$, i.e. if $t\ge \rho$.

\smallskip

In particular, we denote
$\TT=\{\x^{\alpha_1},\ldots,\x^{\alpha_{\d}}\}$, and we assume that $
\TT^*=\{\x^{\alpha_1},\dots,\x^{\alpha_s}\}$,
 the first
$s:=\#(\TT^*)$ elements of $\TT$.

\item $K[{\bf x}]_{t,*}$ denotes the $K$-vector space generated by
$K[\x]_{t}\cup\TT^*$ and $N^*:=\dim(K[\x]_{t,*})$.

\item For $1\le i\le n+1$, $M_{f_i}\in K^{\dim(\RR_i)\times N^*}$ denotes  the transpose of the matrix in the monomial
bases of the composition between the Sylvester multiplication map
and the inclusion $K[\x]_t\to K[\x]_{t,*}$:
$$\begin{array}{cccc}
\mu_{f_i}:& \langle \RR_i\rangle& \to &K[\x]_{t,*} \\
& \x^\alpha&\mapsto&\x^\alpha f_i\end{array}.$$ For later
convenience we order the monomial basis of $K[\x]_{t,*}$ in such a
way that all monomials in $\TT$ precede the monomials in
$K[\x]_{t,*}\smallsetminus\TT$.

 \item $\widetilde M_\SS
\in K^{(N-k)\times (N-k)}$ denotes the Macaulay-Chardin matrix obtained
from
\begin{equation}\label{sylv}
\left[\begin{array}{c}\;\;M_{f_1}\;\;\\[-2mm]\vdots \\[-2mm]\;\;M_{f_{n+1}}\;\;\end{array}\right]
\end{equation}
by deleting the  columns indexed by the monomials in
$\SS\cup\TT^*$.

\item Following \cite{Mac2,Cha}, we define {\em the extraneous
factor}  $\EE(t)$ as the determinant of the square submatrix of
(\ref{sylv}) whose rows are indexed by all those monomials
$\x^\alpha\in {\mathcal R}_i$, $1\leq i\leq n$, such that
$t-d_i-|\alpha|\ge d_{n+1}$ or there exists $j>i$ with
$\alpha_j\geq d_j$, and whose columns are indexed by those
$\x^\alpha$ such that $t-|\alpha|\ge d_{n+1}$ and for some  index
$i$,  $\alpha_i\ge d_i$, or such that there exist at least two
different indexes $1\leq i,j\leq n$ with  $\alpha_i\geq d_i,\
\alpha_j\geq d_j$. It is straightforward to verify that this is
really a square matrix. An important property of $\EE(t)$ is that
it
  neither depends on the coefficients of $f_{n+1}$ nor on $\SS$.

\item $\Delta_{\SS}:=\Delta^{(t)}_{\SS^h}(\fh_1,\dots,\fh_{n+1})$
denotes the {\em order $t$  subresultant of $\fh_1,\dots,\fh_{n+1}$
with respect to}
$\SS^h:=\{\x^{\gamma_1}x_{n+1}^{t-|\gamma_1|},\ldots,\x^{\gamma_k}x_{n+1}^{t-|\gamma_k|}\}.
$  Here, $\fh_i$ denotes the homogenization of $f_i$ by the
 variable $x_{n+1}$. It turns out that by \cite{Cha2} we have
\begin{eqnarray}\label{Extra}
\Delta_{\SS}=\pm\,\frac{|\widetilde M_\SS|}{\EE(t)}.
 \end{eqnarray}

 \item For $1\le i \le n$, $\overline f_i$ is the homogeneous component
 of degree $d_i$ of $f_i$,
  and
 $\overline\Delta_{\TT_j}:=\Delta_{\TT_{j}}^{(j)}(\overline f_1,\dots,\overline f_n)$
 is the order $j$ subresultant of $\overline f_1,\dots,
 \overline f_n$  with respect to $\TT_{j}$.

\item  $\{{\bf \xi}_1,\dots, {\bf \xi}_{\d}\}$ denotes  the
set of all
 common roots of $f_1,\dots, f_n$ in $\overline K^n,$ and
$\VV_{\TT}:=\det({\bf \xi}_j^{\alpha_i})_{1\le i, j\le \d}$, the
generalized  Vandermonde determinant associated to $\TT.$
\end{itemize}

\begin{remark}\label{resutant2} The order $t$ subresultant given in
(\ref{Extra}) generalizes both the univariate case and the usual
multivariate projective resultant as defined for instance in
\cite[Th.~2.3]{CoLiSh2}.
 \begin{enumerate} \vspace{-3mm}\item When $n=1$ and
$t\le d_1+d_2-1$, there are
 no rows and columns of (\ref{sylv})
satisfying the condition that contributes to the extraneous factor
$\EE_t$, and thus  $\EE(t)=1$. Therefore $\Delta_\SS$ of
(\ref{Extra}) coincides with  the univariate order $t$
subresultant of $f$ and $g$ with respect to $\SS$  defined in
Section \ref{univariate}.
\item
 When $t\ge \rho + d_{n+1}$, then $k=0$
since  $\alpha_1<d_1, \dots, \alpha_n<d_n$ imply $|\alpha|\le \rho$,
thus $t-|\alpha|\ge d_{n+1}$.  Therefore $\SS:=\emptyset$. In that case
we recover Macaulay's construction \cite[Th.~p.9 and Th.~4]{Mac2} and
$\Delta_\SS=\pm \Res(f_1^h,\dots, f_{n+1}^h)$.
\end{enumerate}
\end{remark}

We are ready now to state the  main result  of the paper, the multivariate
generalization of Theorem \ref{subresultant1}.

\begin{theorem}\label{subresultant2} Let $f_1, \ldots,
f_{n+1}\in \ K[\bf x]$ and $\{{\bf \xi}_1,\dots, {\bf \xi}_{\d}\}$ be
the set of common roots of $f_1,\dots,f_n$ in $\overline{K}^n$. Then,
under the previous notations, for any $t\in\Z_{\geq0}$ and for any
$\SS=\{\x^{\gamma_1},\ldots,\x^{\gamma_k}\}\subset K[\x]_t$  of
cardinality  $k=\HH_{d_1\dots d_{n+1}}(t)$, the order $t$ subresultant
$\Delta_\SS$  satisfies:
\begin{equation}\label{mff}
\Delta_\SS=\pm \left(\prod_{j=t-d_{n+1}+1}^t \overline
\Delta_{\TT_{j}}\right) \frac{|\OO_\SS|}{\VV_{\TT}},
\end{equation} where
$$\OO_\SS=\left[\begin{array}{ccc}
{\bf \xi}_1^{\gamma_1} &\cdots& {\bf \xi}_{\d}^{\gamma_1}\\[-2mm]
\vdots & & \vdots\\[-2mm]
{\bf \xi}_1^{\gamma_k} &\cdots& {\bf \xi}_{\d}^{\gamma_k}\\[-2mm]
& & \\[-5mm]
\hline
& & \\[-7mm]
{\bf \xi}_1^{\alpha_1} &\cdots& {\bf \xi}_{\d}^{\alpha_1}\\[-2mm]
\vdots & & \vdots \\[-2mm]
{\bf \xi}_1^{\alpha_s} &\cdots& {\bf \xi}_{\d}^{\alpha_s}\\[-2mm]
& & \\[-5mm]
\hline
& & \\[-7mm]
{\bf\xi}_1^{\beta_1}f_{n+1}({\bf\xi}_1)&\cdots&
{\bf\xi}_{\d}^{\beta_1}f_{n+1}({\bf\xi}_{\d})\\[-2mm]
\vdots& &\vdots\\[-2mm]
{\bf \xi}_1^{\beta_r}f_{n+1}({\bf\xi}_1)& \cdots&
{\bf\xi}_{\d}^{\beta_r}f_{n+1}({\bf\xi}_{\d})
\end{array}\right] \ \in \ \overline{K}^{\d\times \d}.$$
\end{theorem}

\begin{proof}
First we check that $\OO_\SS$ is a square matrix, i.e. that $\d=k+s+r$.
This is clear by Formula (\ref{key}) since
 \vspace{-7mm} $$\d-s= \#(
\TT)-\#(\TT^*)=\#(\TT\smallsetminus \TT^*)=\#\{\x^\alpha,|\alpha|\le t,
\alpha_j<d_j \ \forall j\}=k+r.\vspace{-7mm}$$
  In this proof the monomial basis
$\{\x^{\delta_1},\ldots,\x^{\delta_{N^*}}\}$ of $K[\x]_{t,*}$ is
ordered such as was specified in the notations (monomials in  $\TT$
precede the rest of the monomials in $K[\x]_{t,*}$).

 Like in the
univariate case, we define $I_\SS\in K^{(k+s)\times N^*}$ as the
transpose of the matrix of the immersion of the $K$-vector space
generated by $\SS\cup\TT^*$ into $K[\x]_{t,*}$ in the monomial
bases. We set
$$M_\SS:=\left[\begin{array}{ccc}&I_\SS&\\\hline &M_{f_1}&\\ [-2mm]
&\vdots &\\[-2mm]&\;\;M_{f_n}\;\;&\\
\hline&\;\;\;\;\;M_{f_{n+1}}\;\;\;\;\;&\end{array} \right] \in
K^{N^*\times N^*}.$$ ($M_\SS $ is a square matrix by (\ref{MSsquare})
and since $N^*=N+\dim(\TT^*)$.) \par  Furthermore, it is immediate to
verify that $| M_\SS|=\pm\,| \widetilde M_\SS|=\pm \,\EE(t)
\,\Delta_\SS$, where
 $\EE(t)$ denotes the extraneous factor that has been introduced in (\ref{Extra}).

 We set
$$V_{N^*}=\left[\begin{array}{ccc}\xi_1^{\delta_1}  &  \dots &
\xi_{\d}^{\delta_1}\\[-1mm]
\vdots & & \vdots \\[-1mm] \xi_1^{\delta_{N^*}}  &  \dots &
\xi_{\d}^{\delta_{N^*}}
\end{array}\right] \ \in \overline{K}^{N^*\times {\d}} , \ \mbox{and}
\ V_{\d} :=\left[\begin{array}{c|c}&{\bf 0}
\\ \;\;V_{N^*}\;\;&\\[-2mm]
 &\;\;Id \;\;\end{array}\right]
 \in\overline{K}^{N^*\times N^*}$$ and we observe that
$\VV_{\TT}=| V_{\d}|$. We perform the product $M_\SS \,V_{\d}$:
$$ M_\SS\,V_{\d}\,=\,
\begin{tabular}{|c|}
\cline{1-1}
$I_\SS $\\[-1mm]
\cline{1-1}
$M_{f_1}$  \\[-2mm]
\vdots\\[-2mm]
$M_{f_n}$\\[-1mm]
\cline{1-1}
$M_{f_{n+1}}$\\
\cline{1-1}
\end{tabular}
\cdot
\begin{tabular}{|c|c|}
\cline{1-2} &${\bf 0}$ \\[-2mm]$\xi_j^{\delta_i}$
& \\[-2mm]
 &$Id$
\\
\cline{1-2}
\end{tabular}=\begin{tabular}{|c|ccc|}
\cline{1-4}
& & & \\[-8mm]
 $\xi_j^{\gamma_i}$ &  &  $\bf{*}$ & \\
&&&
 \\[-8mm]
 \cline{1-4}
 &&&\\[-8mm]
$\xi_j^{\alpha_i}$ &  &  $\bf{*}$ & \\
& & & \\[-8mm]
\cline{1-4} &&&\\[-8mm]
& &$M'_{f_1}$ & \\[-2mm]
$ {\bf 0}$&
& \vdots&  \\[-2mm]
& &$M'_{f_n}$ &
\\[1mm]
\cline{1-4}
& && \\[-4mm]
 $\xi_j^{\beta_i}f_{n+1}(\xi_j)$ & & ${\bf *}$& \\
& && \\
\cline{1-4}
\end{tabular} \ ,
$$

where $M':=\left[\begin{array}{c}M'_{f_1}\\[-2mm] \vdots \\[-2mm]
\;\;M'_{f_n}\;\;\end{array}\right]$ is the
submatrix of $\left[\begin{array}{c}M_{f_1}\\[-2mm] \vdots \\[-2mm]
\;\; M_{f_n}\;\;\end{array}\right]$ with the same number of rows and
whose columns are
 indexed by all monomials in
 $\x^\alpha \in K[\x]_{t,*}\smallsetminus\TT =
  K[\x]_t\smallsetminus (\TT\smallsetminus \TT^*)=K[\x]_t\smallsetminus \TT$.
 It is immediate to verify that
 $M'$ is a square matrix since, again by  (\ref{key}),
  $\#(\RR_1)+\cdots +\#(\RR_n)=N-k-r=N
 -\#(\TT\smallsetminus\TT^*)=N^*-{\bf d}$.

 We recall that $\#(\TT\smallsetminus \TT^*)=\#\{\x^\alpha,|\alpha|\le t,
\alpha_i<d_i \ \forall i\}$, and  therefore $M'$ is the
Macaulay-Chardin matrix associated to the computation of
$\Delta^{(t)}_{\TT\smallsetminus \TT^*}(f_1,\dots,f_n)$, the order $t$
subresultant of $f_1,\dots,f_n$ with respect to $\TT\smallsetminus
\TT^*$.

To conclude the proof we are left to prove that
$$|M'|=\pm \EE(t) \,\left(\prod_{j=t-d_{n+1}+1}^t
\overline \Delta_{\TT_{j}}\right).$$

This was proven in  \cite[p.14]{Mac2} (see also the proof of
\cite[Lem.~1]{Cha2} and \cite[Thm.~5.2]{DaJe2004}). For the
reader's convenience, we rewrite the proof here.


We reorganize the matrix $M'$ as follows: we recall that the columns
correspond to monomials $\x^\alpha \in K[\x]_{t}\smallsetminus\TT$ and we
index the
 columns by graded descending order, first all monomials of degree
 $t$ in $K[\x]_{t}\smallsetminus\TT$,
 then
all monomials of degree $t-1$ in $K[\x]_{t}\smallsetminus\TT$, and so on,
up to all monomials of degree $t-d_{n+1}+1$. Finally, we put in the
last block all monomials of degree bounded by $t-d_{n+1}$. The rows
correspond to $\RR_i$ for $1\le i\le n$. We also index them by
graded descending order: first all monomials of degree $t-d_i$ in
$\RR_i$ for $1\le i\le n$, then all monomials of degree $t-d_i-1$ in
$\RR_i$, $1\le i\le n$, and so on up to all monomials of degree
$t-d_i-d_{n+1}+1$ in $\RR_i$, $1\le i\le n$. In the last block we
put all monomials of degree bounded by $t-d_i-d_{n+1}$ in $\RR_i$,
$1\le i\le n$.

With this ordering $M'$ has
 a block structure:
\begin{equation}\label{cito}
M'=\left[\begin{array}{cccc}M_t & {\bf *} & {\bf *} & {\bf *}
\\ &  \ddots&   {\bf *} & {\bf *}\\
& &  M_{t-d_{n+1}+1}& {\bf *}\\
{\bf{0}} & & & E
\end{array}\right],
\end{equation}
where the square matrix $M_j$ corresponds to the coefficients of the
terms of degree $j$ of $\x^\alpha f_i,$ where $|\alpha|=j-d_i$, that
is, the coefficients of $\x^\alpha \overline{f}_i$ except those
corresponding to terms in $\TT_j.$

Hence $M_j$ is the Macaulay-Chardin matrix associated to the
$j$-subresultant $\overline{\Delta}_{\TT_j}$ of
$\overline{f}_1,\dots, \overline{f}_n$ with respect to $\TT_j$
(\cite{Cha}) and it turns out that
$$| M_j |=\EE_j\,\overline{\Delta}_{\TT_j},$$ where $\EE_j$ is the extraneous factor
associated to this construction, that we recall  only depends on $j$
and not on the set $\TT_j.$

But it turns out that the extraneous factor $\EE(t)$
 has a block structure similar to (\ref{cito}) (see \cite{Mac2,Cha2,DaJe2004}).
We have, with our notation:
\begin{equation}\label{rar}
\EE(t)=|E|\,\prod_{j=t-d_{n+1}+1}^t \EE_j,
\end{equation}
(see \cite[Th.~6]{Mac2}). This concludes the proof of the Theorem.
\end{proof}

\smallskip
\begin{remark}\label{kokk}
The reason why we cannot allow $\TT_j$ to be \textit{any} subset of
monomials of degree $j$ for  $j\le t-d_{n+1}+1$ is the factorization
formula on the right hand side of (\ref{rar}), where the $\EE_j$'s
involved in the product are only those corresponding to $j$ satisfying
$t-d_{n+1}+1\leq j\leq t.$
 This is not just a technical obstruction.
If we could pick any $\TT_j$ for every $j,$ then setting
$t:=\rho+d_{n+1}$,  the Poisson formula for the resultant
$\Res(\fh_1, \ldots, \fh_{n+1})$  would read as follows
$$
\frac{\left|\begin{array}{ccc}
{\bf\xi}_1^{\beta_1}&\cdots& {\bf\xi}_{\d}^{\beta_1}\\
\vdots& &\vdots\\
{\bf \xi}_1^{\beta_r}& \cdots&{\bf\xi}_{\d}^{\beta_r}
\end{array}\right|}{V_\TT}
\Res(\overline f_1,\dots,\overline
 f_n)^{d_{n+1}}\prod_{\xi\in V_{{\overline K}^n}(f_1, \ldots, f_n)}
f_{n+1}(\xi),$$ which is obviously false in general since the
fraction does not cancel unless $\TT=\RR_{n+1},$ i.e. $\TT_j$ is
defined as in (\ref{prof}).
\end{remark}

\smallskip

Like in the univariate case, we illustrate Theorem
\ref{subresultant2} with a specific example.

\begin{example}
Let $n=2$, $d_1=d_2=d_3=2$ and $t=t^*=2$.

 Here $k=\#\{x_1,x_2,x_1x_2\}=3$, $\RR_1=\RR_2=\RR_3=\{1\}$ and
 $\TT=\{1,x_1,x_2,x_1x_2\}$.

We fix the ordered monomial basis $(1,x_1,x_2,x_1x_2,x_1^2,x_2^2)$
of $K[\x]_2$ and
 \begin{eqnarray*}
{f}_{{1}}\, &=& \,a_{{0}}+a_{{1}}{x_1}+a_2\,{x_2}+
a_{{3}}{x_1x_2}+a_{{4}}x_1^{2}+a_{{5}}{x_2^2}\\
{f}_{{2}}\, &=& \,b_{{0}}+b_{{1}}{x_1}+b_2\,{x_2}+
b_{{3}}{x_1x_2}+b_{{4}}x_1^{2}+b_{{5}}{x_2^2}\\
{f}_{{3}}\, &=& \,c_{{0}}+c_{{1}}{x_1}+c_2\,{x_2}+
c_{{3}}{x_1x_2}+c_{{4}}x_1^{2}+c_{{5}}{x_2^2}.
 \end{eqnarray*}

  Then
 $$
 \left[\begin{array}{c}M_{f_1} \\[-2mm]M_{f_2}\\[-2mm]M_{f_{3}}\end{array}\right]= \left[ \begin{array}{cccccc}
a_{{0}}&a_{{1}}&a_{{2}}&a_{{3}}&a_{{4}}&a_{{5}}\\[-2mm]
b_{{0}}&b_{{1}}&b_{{2}}&b_{{3}}&b_{{4}}&b_{{5}}\\[-2mm]
c_{{0}}&c_{{1}}&c_{{2}}&c_{{3}}&c_{{4}}&c_{{5}}\end{array} \right] .
$$
We choose $\SS:=\{x_1,x_1x_2,x_1^2\}$. Then
$$
\Delta_\SS=c_{{0}}(a_{{2}}b_{{5}}-a_{{5}}b_{{2}})-c_2(a_{{0}}b_{{5}}-a_{{5}}b_{{0}})+
c_{{5}}(a_{{0}}b_{{2}}-a_{{2}}b_{{0}}) .
$$
On the other hand, if $V_{\overline
K}(f_1,f_2)=\{\xi_1,\xi_2,\xi_3,\xi_4\} $ with
$\xi_j=(\xi_{j1},\xi_{j2})$ for $1\le j\le 4$, then
 $$
 \OO_\SS\, = \, \left[ \begin {array}{cccc}
\xi_{{11}}&\xi_{{21}}&\xi_{{31}}&\xi_{{41}}\\[-2mm]
\noalign{\medskip}\xi_{{11}}\xi_{{12}}&\xi_{{21}}\xi_{{22}}&\xi_{{31}}
\xi_{{32}}&\xi_{{41}}\xi_{{42}}\\[-2mm]
\noalign{\medskip}\xi_{{11}}^{2}&\xi_{{21}}^{2}&\xi_{{31}}^{2}&\xi_{{41}}^{2}\\[-2mm]
\noalign{\medskip}f_3(\xi_1)&f_3(\xi_2)&f_3(\xi_3)&f_3(\xi_4)
\end
{array} \right] .
$$
Therefore, if we set $V$ for the generalized Vandermonde matrix on
$\xi_1,\dots,\xi_4$ corresponding to the sequence of monomials
$1,x_1, x_2, x_1x_2,x_1^2, x_2^2$, i.e.
$$
V:= \, \left[ \begin {array}{cccc} 1 & 1& 1& 1\\[-3mm]
\xi_{{11}}&\xi_{{21}}&\xi_{{31}}&\xi_{{41}}\\[-2mm]
\xi_{{12}} &\xi_{{22}}&\xi_{{32}} &\xi_{{42}}\\[-2mm]
\noalign{\smallskip} \xi_{11}\xi_{{12}}&\xi_{{21}}\xi_{{22}}&\xi_{{31}}
\xi_{{32}}&\xi_{{41}}\xi_{{42}}\\[-2mm]
\noalign{\medskip}\xi_{{11}}^{2}&\xi_{{21}}^{2}&\xi_{{31}}^{2}&\xi_{{41}}^{2}\\[-2mm]
\noalign{\medskip}\xi_{{12}}^{2}&\xi_{{22}}^{2}&\xi_{{321}}^{2}&\xi_{{42}}^{2}
\end
{array} \right] \in \overline{K}^{6\times 4},
$$
and $V_{i,j}$, $0\leq i<j\leq 5$, for  the square submatrix
 obtained from $V$  deleting the  $i$-th and $j$-th rows
(we adopt the convention of numbering the rows  from $0$ to $5$ like
the coefficients of the $f_i$'s), we conclude that
$$
|\OO_\SS|=-c_0\,|V_{2,5}|+c_2\,|V_{0,5}|+c_5\,|V_{0,2}|.
$$
 Also,  with this
notation $V_{4,5}$ is the  Vandermonde matrix corresponding to
$\TT$.

Now, since the only non-trivial homogeneous subresultant $\overline
\Delta_{\TT_{j}}$ in (\ref{mff}) is for $\TT_2=\{x_1x_2\}$, and is
equal to
 $$\overline \Delta_{\TT_{2}}=a_4b_5-a_5b_4,$$
 Theorem \ref{subresultant2} states that
 \begin{eqnarray*}
c_{{0}}(a_{{2}}b_{{5}}-a_{{5}}b_{{2}})-c_2(a_{{0}}b_{{5}}-a_{{5}}b_{{0}})+
c_{{5}}(a_{{0}}b_{{2}}-a_{{2}}b_{{0}})\\= \pm (a_4b_5-a_5b_4)
\left(-c_0\frac{|V_{2,5}|}{|V_{4,5}|}+c_2\frac{|V_{0,5}|}{|V_{4,5}|}+
c_5\frac{|V_{0,2}|} {|V_{4,5}|}\right).
\end{eqnarray*}
Indeed, we show below that this equality  holds  since   for any
$i<j$ and $k<l$:
\begin{equation}\label{exampleeq}
(-1)^{i+j}\frac{|V_{i,j}|}{a_{{i}}b_{{j}}-a_{{j}}b_{{i}}}=(-1)^{k+l}\frac{|V_{k,l}|}{a_kb_l-a_lb_k}
.\end{equation}

If for $0\le i,j\le 5$, we set $I_{i,j}\in K^{4\times 6}$ a
$4$-identity matrix with added $0$ columns for column $i$ and column
$j$, and $I^{i,j}\in K^{6\times 2}$ the matrix with $4$ null rows and
the identity matrix plugged in rows $i$ and $j$, we observe that

\begin{eqnarray*}
\begin{tabular}{|ccc|}
\cline{1-3} & & \\[-7mm]
&$I_{i,j}$&\\
& & \\[-8mm]
\cline{1-3}&&\\[-8mm]
&$M_{f_1}$ & \\[-2mm]
&$M_{f_2}$&\\
 \cline{1-3}
\end{tabular}
\cdot
\begin{tabular}{|ccc|c|}
\cline{1-4} & && \\[-2mm]  & $V$ & &${I^{k,l}}$    \\[-2mm]  & && \\
\cline{1-4}
\end{tabular}
&=&
 \begin{tabular}{|ccc|ccc|}
\cline{1-6} & & && & \\[-7mm]
   & $V_{i,j}$& &&$*$&\\
&& & &&\\[-8mm]
\cline{1-6}&&&&&\\[-8mm]
 & &&$a_k$& &$a_l$\\[-2mm]
&{\bf 0}&&$b_k$&&$b_l$\\[-8mm]& & && & \\
 \cline{1-6}
\end{tabular}
\end{eqnarray*}

since $f_1(\xi_j)=f_2(\xi_j)=0$, $1\le j\le 4$.  Thus, taking
determinants on both sides,
$$(-1)^{5-j+4-i}(a_ib_j-a_jb_i)\cdot(-1)^{k+l-1}|V_{k,l}|=|V_{i,j}|\cdot (a_kb_l-a_lb_k),$$
and we obtain  (\ref{exampleeq}).

Applying this to our case, we  conclude that here
$$\Delta_\SS=-\,
\left(\prod_{j=t-d_{n+1}+1}^t \overline \Delta_{\TT_{j}}\right)
\frac{|\OO_\SS|}{\VV_{\TT}}.
$$
{\hfill\mbox{$\Box$}}
\end{example}

\smallskip

Next, we   recover
 Theorem \ref{subresultant1} in the univariate case:

\begin{observation}
For $n=1$, by setting $f_1:=g$ and $f_2:=f$, as
$\overline{f_1}=b_{d_2}x^{d_2},$ it turns out that
$$\overline{\Delta}_{\TT_j}=\left\{\begin{array}{ccc}
b_{d_2}& \mbox{if} & j\geq d_2 \\
1 & \mbox{if} & j<d_2.
\end{array}\right.$$
So, if $t\geq d_2,$ then $\prod_{j=t-d_1+1}^t \overline
\Delta_{\TT_{j}}={b_{d_2}}^{t-d_2+1}.$ If $t<d_2,$ the product of
subresultants equals $1 $.  \hfill\mbox{$\Box$}
\end{observation}

\smallskip In the particular case $t=\rho+d_{n+1}$, Theorem
\ref{subresultant2} gives a new proof  for  the Poisson product
formula for the multivariate resultant (see \cite{CoLiSh2}):


\begin{corollary} \ {\em (Poisson product formula)}
$$ \Res(\fh_1, \ldots, \fh_{n+1})=\pm\Res(\overline f_1,\dots,\overline
 f_n)^{d_{n+1}}\prod_{\xi\in V_{{\overline K}^n}(f_1, \ldots, f_n)}
f_{n+1}(\xi).$$
\end{corollary}

\begin{proof}
We apply Remark \ref{resutant2} (2) for $t:=\rho + d_{n+1}$ to Theorem
\ref{subresultant2}.
We observe that by the same remark,  for $j> \rho$, i.e. for $j\ge t-d_n+1$,
$\overline{\Delta}_{\TT_j}=\Res(\overline f_1,\dots,\overline
 f_n)$. We conclude  that
$\OO_\SS$ equals $\left(\prod_{\xi\in V_{{\overline K}^n}(f_1,
\ldots, f_n)} f_{n+1}(\xi)\right)$ times the generalized Vandermonde
matrix whose determinant equals $\VV_{\TT}$.
\end{proof}

\smallskip

We end this paper by giving the multivariate version of Corollary
\ref{koko}, i.e. a discrete Wr\'onskian type expression for the {\em
generalized subresultant polynomial}:

\begin{equation}
\label{chardinpol2} s(\x):=\sum_{j=0}^k \Delta_{\SS_{j}}\x^{\gamma_j},
\end{equation}
  defined for a fixed $t\in \N$ and  $k:=\HH_{d_1\dots d_{n+1}}(t)$,
  under the usual notations, \par
  $\SS:=\{\x^{\gamma_j},0\le j\le k\}\subset K[\x]_t$ and
$\SS_{j}:=\SS \smallsetminus \{ \x^{\gamma_j}\}$.

It turns out that
 $s(\x)$ belongs to the ideal generated by the $f_i$'s (see
 \cite{Cha}), and
the following result can be proved mutatis mutandis the proof of
Corollary \ref{hong}.

\begin{corollary} Let $f_1, \ldots, f_{n+1}\in K[\x]$
and  $s(\x)$ be the generalized subresultant polynomial defined in
(\ref{chardinpol2}). Then,
   we have
{\small
$$s(\x)=\pm {\VV_\TT}^{-1} \left(\prod_{j=t-d_{n+1}+1}^t \overline
\Delta_{\TT_{j}}\right) \left|\begin{array}{cccc}
\x^{\gamma_0}& {\bf \xi}_1^{\gamma_0}& \dots & {\bf \xi}_{\d}^{\gamma_0}\\[-2mm]
\x^{\gamma_1} &{\bf\xi}_1^{\gamma_1}& \dots & {\bf \xi}_{\d}^{\gamma_1}\\[-2mm]
\vdots& \vdots &  & \vdots\\[-2mm]
\x^{\gamma_k}& {\bf\xi}_1^{\gamma_k} &\dots & {\bf\xi}_{\d}^{\gamma_k}\\[-2mm]
0&{\bf\xi}_1^{\xi_1}f_{n+1}({\bf\xi}_1)&\cdots& {\bf\xi}_{\d}^{\xi_1}f_{n+1}
({\bf\xi}_{\d})\\[-2mm]
\vdots&\vdots& &\vdots\\[-2mm]
0& {\bf\xi}_1^{\xi_r}f_{n+1}({\bf\xi}_1)& \cdots&
{\bf\xi}_{\d}^{\xi_r}f_{n+1}({\bf\xi}_{\d})
\end{array}\right|.$$}
\hfill\mbox{$\Box$}
\end{corollary}

\begin{remark}
If $\gcd(\SS)\in\SS,$ then one can reduce the previous determinant,
as in Corollary \ref{koko}.
\end{remark}


\end{document}